\documentclass[12pt]{article}
\usepackage{graphicx}
\usepackage{amsmath}
\usepackage{epsf}
\usepackage{graphpap}

\numberwithin{equation}{section}

\def\be{\begin{equation}}
\def\ee{\end{equation}}

\date{February 02, 2008}
\begin{document}

\title{Use of Complex Lie Symmetries for Linearization of Systems of \\
Differential Equations -- I:\\
Ordinary Differential Equations}
\author{S. Ali$^{a}$, F. M. Mahomed$^{b}$, A. Qadir$^{a}$ \\
%EndAName
$^{a}${\small Center For Advanced Mathematics and Physics, }\\
{\small National University of Sciences and Technology,} \\
{\small Rawalpindi, Pakistan.}\\
$^{b}${\small Centre for Differential Equations, Continuum }\\
{\small Mechanics and Applications, School of Computational }\\
{\small and Applied Mathematics, University of} {\small the}\\
{\small Witwatersrand, South Africa.}\\
}
\maketitle

\begin{abstract}
The Lie linearizability criteria are extended to complex functions for
ordinary differential equations. The linearizability of complex ordinary
differential equations is used to study the linearizability of corresponding
systems of two ordinary differential equations. The transformations that map
a system of two nonlinear ordinary differential equations into systems of
linear ordinary differential equations are obtained from complex
transformations. Invariant criteria for linearization are given for second
order complex differential equations in terms of the coefficients of the
equations, as well as the corresponding real system, which provide
procedures for writing down the solutions of the equations. Illustrative
examples are given and discussed.
\end{abstract}

\section{\textbf{Introduction}}

\qquad Lie provided linearizability criteria for scalar second order real
ordinary differential equations (RODEs) \cite{lie}. The criteria help in
constructing those transformations that map a nonlinear second order RODE
into a linear second order RODE. He proved that to be linearizable, a second
order RODE must be at most cubic in the first derivative together with the
condition that the coefficients satisfy an overdetermined system of
equations (see e.g., $[6-8]$, [$10-12$], [15], [20] \ [21]). Tress\'{e} \cite%
{tre} gave necessary and sufficient conditions for linearization of RODEs by
calculating the two relative invariants of the equivalence group of
transformations, which are equivalent to the Lie conditions \cite{lea}. The
Lie approach has been extended in many ways. The linearization of higher
order RODEs and systems of RODEs has been discussed by various authors (see
e.g., [9], $[13-18])$.

In a companion paper, henceforth called P-II\ \cite{Ali2}, the Lie approach
of linearization of RODEs is extended to linearization of complex ordinary
differential equations (CODEs) via complex point transformations. A CODE is
a differential equation in a complex function of a single complex variable.
Here we are interested in developing the insights obtained by considering
ODEs for complex functions of a single real variable (which we call
r-CODEs). Whereas linear r-CODEs are trivial, nonlinear r-CODEs couple the
two components of the complex function to work what Penrose calls
\textquotedblleft \textit{complex magic}\textquotedblright\ \cite{Pen}. An
r-CODE can be linearized in the same way as a RODE via invertible
complex transformations. It yields
two RODEs. The linearization of such systems of RODEs follows directly from
the linearization of the corresponding r-CODEs after decomposing the
complex transformation into two real point transformations.

The paper is sectioned as follows. In the second section equivalent r-CODEs
are discussed i.e. those which can be transformed into one another via
invertible mappings. The third section includes the derivation of Lie
compatibility conditions and some other results for linearization of
r-CODEs. Illustrative examples are also discussed in the same section.
Finally, we give a brief summary and discussion of our complex Lie symmetry
(CLS) approach.

\section{\textbf{Equivalent r-CODEs}}

\qquad Two r-CODEs are (locally) equivalent via an invertible complex
transformation if one can be transformed into the other by an invertible
complex transformation. An r-CODE is a system of two RODEs in two
unknown functions of one variable. If an r-CODE is equivalent to another
r-CODE via invertible complex transformations then the system of RODEs
corresponding to that r-CODE is also equivalent to the other system of
RODEs. The invertible complex transformation yields two real invertible
transformations which can be used to transform one system of RODEs into
another system. We illustrate this fact in the following examples. We shall
write $u(x)$ as a complex function $u$ of a single real variable, $x.$

Every first-order r-CODE of the form 
\begin{equation}
u^{\prime }(x)=w(x,u),  \tag{1}
\end{equation}%
is equivalent to a simple r-CODE, $u^{\prime }=0,$ via a complex
transformation which is the same as the case of RODEs \cite{mah1}. Further,
every\textbf{\ }\textit{linear }second-order r-CODE is also equivalent to
its simple analogue $u^{\prime \prime }(x)=0,$ via some complex
transformation. For example, the complexified Ricatti equation 
\begin{equation}
u^{\prime }(x)+u^{2}=0,  \tag{2}
\end{equation}%
is transformable to $U^{\prime }=0,$ by means of 
\begin{equation}
\chi =x,~U=(1/u)-x.  \tag{3}
\end{equation}%
Similarly, the complexified simple harmonic oscillator equation 
\begin{equation}
u^{\prime \prime }(x)+u=0,  \tag{4}
\end{equation}%
can be transformed into $U^{\prime \prime }=0,$ via an invertible
complex transformation%
\begin{equation}
\chi =\tan x,~U=u\sec x  \tag{5}
\end{equation}%
in the domain $(-\pi /2,\pi /2).$

Now, we decompose the above r-CODEs into systems of RODEs corresponding to
them by 
\begin{equation}
u(x)=y(x)+iz(x).  \tag{6}
\end{equation}%
\newline
The system of RODEs corresponding to a first-order r-CODE is 
\begin{equation}
y^{\prime }=w_{1}(x,y,z),\text{ \ }z^{\prime }=w_{2}(x,y,z).  \tag{7}
\end{equation}%
This system is transformable into%
\begin{equation}
\Upsilon ^{\prime }=0,\text{ \ \ \ }\zeta ^{\prime }=0,  \tag{8}
\end{equation}%
by the invertible transformation derived from the complex
transformation.\textbf{\ }Similarly, a general linear second-order r-CODE
yields the following system of RODEs 
\begin{equation}
y^{\prime \prime }=w_{1}(x,y,z,y^{\prime },z^{\prime }),  \tag{9}
\end{equation}%
\begin{equation}
z^{\prime \prime }=w_{2}(x,y,z,y^{\prime },z^{\prime }),  \tag{10}
\end{equation}%
where both $w_{1}$ and $w_{2}$ are two such real functions that do not give
rise to a nonlinear system. This system can be transformed into%
\begin{equation}
\Upsilon ^{\prime \prime }=0,\text{ \ \ }\zeta ^{\prime \prime }=0.  \tag{11}
\end{equation}%
\quad Notice that the general solution of (4) is 
\begin{equation}
u(x)=\alpha \cos x+\beta \sin x,  \tag{12}
\end{equation}%
where $\alpha$ and$~\beta$ are complex constants such that $\alpha =\alpha _{1}+i\alpha _{2}$ and $~\beta =\beta _{1}+i\beta
_{2}$ despite the fact that $x$ is a real variable, because the
solution is obtained from complex integration. By using (6), the above
solution becomes%
\begin{equation}
y(x)=\alpha _{1}\cos x+\beta _{1}\sin x,~z(x)=\alpha _{2}\cos x+\beta
_{2}\sin x,  \tag{13}
\end{equation}%
which is in fact the solution of the system%
\begin{equation}
y^{\prime \prime }=-y,~z^{\prime \prime }=-z,  \tag{14}
\end{equation}%
corresponding to r-CODE (4) indicating the consistency of r-CODEs and their
solutions. It is important to note that the solutions of r-CODEs yield
solutions of systems of ODEs corresponding to them.
\newline
\newline
\textbf{Examples}\newline
\textbf{1.} The complexified Ricatti equation (2) yields two RODEs%
\begin{equation}
y^{\prime }=-y^{2}+z^{2},\text{ \ }z^{\prime }=-2yz,  \tag{15}
\end{equation}%
which is a special case in two dimensions of the generalized Ricatti system.
Note that what seems to be a trivial complex transformation yields the
following non-trivial real transformation%
\begin{equation}
\chi =x,\text{~}\Upsilon =\frac{y}{y^{2}+z^{2}}-x,~\zeta =\frac{-z}{%
y^{2}+z^{2}},  \tag{16}
\end{equation}%
that maps system (15) into%
\begin{equation}
\Upsilon ^{\prime }=0,\text{ \ }\zeta ^{\prime }=0.  \tag{17}
\end{equation}%
\textbf{2.} The decoupled linear harmonic oscillator%
\begin{equation}
y^{\prime \prime }=-y,\text{ \ }z^{\prime \prime }=-z,  \tag{18}
\end{equation}%
can be transformed into%
\begin{equation}
\Upsilon ^{\prime \prime }=0,\text{ \ }\zeta ^{\prime \prime }=0,  \tag{19}
\end{equation}%
via the real transformation%
\begin{equation}
\chi =\tan x,\text{ \ }\Upsilon =y\sec x,\text{ \ }\zeta =z\sec x.  \tag{20}
\end{equation}%
This transformation could have been easily guessed. However, the
\textquotedblleft complex magic\textquotedblright\ is that it can be \textit{%
derived} by the use of complex functions of a real variable.

Thus, one can use simple complex transformations to map systems of two
RODEs into systems of simple RODEs which we state as equivalent systems of
RODEs corresponding to some r-CODE.

\section{\textbf{Lie Conditions for Complex Functions}}

\qquad Lie proved that a second order scalar RODE, which is at most cubic in
its first derivative, is linearizable provided it satisfies four
differential constraints involving two auxiliary variables that Tress\`{e}
reduced to two differential conditions \cite{lie, tre}. We extend the Lie
theorem to r-CODEs by considering an r-CODE of the form%
\begin{equation}
u^{\prime \prime }(x)=A(x,u)u^{\prime 3}+B(x,u)u^{\prime 2}+C(x,u)u^{\prime
}+D(x,u),  \tag{21}
\end{equation}%
where $A,~B,~C$ and $D$ are complex valued functions.

Now, for completness, we restate several known results \cite{mah1} extended
to r-CODEs as Theorem 1. \newline
\newline
\textbf{Theorem 1. } The following statements are equivalent.\newline
\textbf{1.} A scalar second order r-CODE is linearizable via a complex
transformation;\newline
\textbf{2.} Equation (21) has a maximal $8-$dimensional complex Lie algebra;%
\newline
\textbf{3.} The Tress\'{e} relative invariants \newline
\begin{gather}
I_{1}=w_{u^{\prime }u^{\prime }u^{\prime }u^{\prime }},  \tag{22} \\
I_{2}=\frac{d^{2}}{dx^{2}}w_{u^{\prime }u^{\prime }}-4\frac{d}{dx}%
w_{u^{\prime }u}-3w_{u}w_{u^{\prime }u^{\prime }}+6w_{u~u}+w_{u^{\prime
}}(4w_{u^{\prime }u}-\frac{d}{dx}w_{u^{\prime }u^{\prime }}),  \tag{23}
\end{gather}%
both vanish identically for (21);\newline
\textbf{4.} The coefficients $A$ to $D$ in (21) satisfy the the Lie type
condition in the complex domain%
\begin{equation}
k_{x}=kK-AD+\frac{1}{3}C_{u}-\frac{2}{3}B_{x},  \tag{24}
\end{equation}%
\begin{equation}
k_{u}=-k^{2}-Bk-AK-A_{x}-AC,  \tag{25}
\end{equation}%
\begin{equation}
K_{x}=K^{2}+Dk+CW-D_{u}+BD,  \tag{26}
\end{equation}%
\begin{equation}
K_{u}=-Kk+AD+\frac{1}{3}B_{x}-\frac{2}{3}C_{u},  \tag{27}
\end{equation}%
where $k$ and $K$ are auxiliary complex functions;\newline
\textbf{5.} The coefficients in (21) namely $A$ to $D$ also satisfy the Lie
compatibility conditions%
\begin{equation}
3A_{xx}+3A_{x}C-3A_{u}D+C_{uu}-6AD_{u}+BC_{u}-2BB_{x}-2B_{xu}=0,  \tag{28}
\end{equation}%
\begin{equation}
6A_{x}D-3B_{u}D+3AD_{x}+B_{xx}-2C_{xu}-3BD_{u}+3D_{{u}u}+2CC_{u}-CB_{x}=0; 
\tag{29}
\end{equation}%
\textbf{6.} Equation (21) has two commuting symmetries $\mathbf{Z}_{1},~%
\mathbf{Z}_{2},$ with$~\mathbf{Z}_{1}=\rho (x,u)\mathbf{Z}_{2},$ for a
nonconstant complex function $\rho ,$ such that a point transformation $%
\zeta =\zeta (x,u),U=U(x,u),$ which brings $\mathbf{Z}_{1}$ and $\mathbf{Z}%
_{2}$ to their canonical form%
\begin{equation}
\mathbf{Z}_{1}=\frac{\partial }{\partial U},~\mathbf{Z}_{2}=\zeta \frac{%
\partial }{\partial U},  \tag{30}
\end{equation}%
reduces the equation to the linear form $U^{\prime \prime }=W(\zeta );$%
\newline
\textbf{7.} Equation (21) has two noncommuting symmetries $\mathbf{Z}_{1},~%
\mathbf{Z}_{2},$ in a suitable basis with$~[\mathbf{Z}_{1},~\mathbf{Z}_{2}]=%
\mathbf{Z}_{1},~\mathbf{Z}_{1}=\rho (x,u)\mathbf{Z}_{2},$ for a nonconstant
complex function $\rho ,$ such that a point change of variables $\zeta
=\zeta (x,u),U=U(x,u),$ which brings $\mathbf{Z}_{1}$ and $\mathbf{Z}_{2}$
to their canonical form%
\begin{equation}
\mathbf{Z}_{1}=\frac{\partial }{\partial U},~\mathbf{Z}_{2}=U\frac{\partial 
}{\partial U},  \tag{31}
\end{equation}%
reduces the equation to the linear form $U^{\prime \prime }=U^{\prime
}W(\zeta );$\newline
\textbf{8.} Equation (21) has two commuting symmetries $\mathbf{Z}_{1},~%
\mathbf{Z}_{2},$ with$~\mathbf{Z}_{1}\neq \rho (x,u)\mathbf{Z}_{2},$ for a
nonconstant complex function $\rho ,$ such that a point transformation $%
\zeta =\zeta (x,u),U=U(x,u),$ which brings $\mathbf{Z}_{1}$ and $\mathbf{Z}%
_{2}$ to their canonical form%
\begin{equation}
\mathbf{Z}_{1}=\frac{\partial }{\partial \zeta },~\mathbf{Z}_{2}=\frac{%
\partial }{\partial U},  \tag{32}
\end{equation}%
reduces the equation to one which is at most cubic in the first derivative;%
\newline
\textbf{9.} Equation (121) has two noncommuting symmetries $\mathbf{Z}_{1}$
and $~\mathbf{Z}_{2}$ in a suitable basis with$~[\mathbf{Z}_{1},~\mathbf{Z}%
_{2}]=\mathbf{Z}_{1}$with$~\mathbf{Z}_{1}\neq \rho (x,u)\mathbf{Z}_{2},$ for
a nonconstant complex function $\rho ,$ such that a point transformation $%
\zeta =\zeta (x,u),U=U(x,u),$ which brings $\mathbf{Z}_{1}$ and $\mathbf{Z}%
_{2}$ to their canonical form%
\begin{equation}
\mathbf{Z}_{1}=\frac{\partial }{\partial U},~\mathbf{Z}_{2}=\zeta \frac{%
\partial }{\partial \zeta }+U\frac{\partial }{\partial U},  \tag{33}
\end{equation}%
reduces the equation to%
\begin{equation}
\zeta U^{\prime \prime }=aU^{^{\prime }3}+bU^{\prime 2}+(1+\frac{b^{2}}{3a}%
)U^{\prime }+\frac{b}{3a}+\frac{b^{3}}{27a^{2}},  \tag{34}
\end{equation}%
where $a$ $(\neq 0)$ and $b$ are complex constants.

\section{\textbf{Lie Conditions for Systems of Two RODEs }}

\qquad The fact that a second order r-CODE is linearizable if it is at most
cubic in its first derivative, i.e. (28) and (29) both hold, implies that a
system of two RODEs corresponding to that r-CODE, which is at most cubic in
the first derivative is also linearizable with the relevant conditions on
the coefficients. We can easily obtain the real transformations by
decomposing each of the complex transformation that linearize an
r-CODE. The general form of a system of RODEs corresponding to (21) is given
by 
\begin{equation}
y^{\prime \prime }=A_{1}(y^{\prime 3}-3y^{\prime }z^{\prime
2})-A_{2}(3y^{\prime 2}z^{\prime }-z^{\prime 3})+B_{1}(y^{\prime
2}-z^{\prime 2})-2B_{2}y^{\prime }z^{\prime }+C_{1}y^{\prime }-C_{2}\zeta
^{\prime }+D_{1},  \tag{35}
\end{equation}%
\begin{equation}
z^{\prime \prime }=A_{1}(3y^{\prime 2}z^{\prime }-z^{\prime
3})+A_{2}(y^{\prime 3}-3y^{\prime }z^{\prime 2})+2B_{1}y^{\prime }z^{\prime
}+B_{2}(y^{\prime 2}-z^{\prime 2})+C_{2}y^{\prime }+C_{1}\zeta ^{\prime
}+D_{2},  \tag{36}
\end{equation}%
where the coefficents $A_{i},B_{i,}C_{i}$ and $D_{i}$ are function of $x,y,z.
$\newline
We obtain results for systems of RODEs by decomposing corresponding r-CODEs
and state them as Theorem 2 with the following notation%
\begin{equation}
\mathbf{Z}_{i} =\mathbf{X}_{i}+\iota \mathbf{Y}_{i}\text{, \ \ \ \ \ (}%
i=1,2)  \tag{37} \\
\end{equation}%
\begin{equation}
\rho (x) =\rho _{1}(x)+\iota \rho _{2}(x),  \tag{38}
\end{equation}%
where $\rho _{1}(x)$ and $\rho _{2}(x)$ are nonconstant real functions and
the complex transformation $(x,u)\longrightarrow (\chi ,U),$ being
equivalent to the real transformation $(x,y,z)\longrightarrow (\chi
,\Upsilon ,\zeta ).$\newline
\newline
\textbf{Theorem 2. }The following statements are equivalent.\newline
\textbf{1.} The system of RODEs (35) \ (36) is linearizable via real
transformations;\newline
\textbf{2.} The coefficients in (35) and (36) satisfy%
\begin{gather}
3A_{xx}^{1}+3C^{1}A_{x}^{1}-3A_{x}^{2}C^{2}-3A_{y}^{1}D^{1}-3D^{1}A_{z}^{2}+3D^{2}A_{y}^{2}-3D^{2}A_{z}^{1}+3A^{1}C_{x}^{1}+
\notag \\
3A^{2}C_{x}^{1}+C_{yy}^{1}-C_{zz}^{1}+2C_{yz}^{2}-6A^{1}D_{y}^{1}-6A^{1}D_{z}^{2}+6A^{2}D_{y}^{2}-6A^{2}D_{z}^{1}+
\notag \\
B^{1}C_{y}^{1}+B^{1}C_{z}^{2}-B^{2}C_{y}^{2}+B^{2}C_{z}^{1}-2B^{1}B_{x}^{1}+2B^{2}B_{x}^{2}-2B_{xy}^{1}-2B_{xz}^{2}=0,
\tag{39}
\end{gather}%
\begin{gather}
3A_{xx}^{2}+3C^{2}A_{x}^{1}+3A_{x}^{2}C^{1}-3D^{2}A_{y}^{1}-3D^{2}A_{z}^{2}-3D^{1}A_{y}^{2}+3D^{1}A_{z}^{1}+3A^{2}C_{x}^{1}-
\notag \\
3A^{1}C_{x}^{1}+C_{yy}^{2}-C_{zz}^{2}-2C_{yz}^{1}-6A^{2}D_{y}^{1}-6A^{2}D_{z}^{2}-6A^{1}D_{y}^{2}+6A^{1}D_{z}^{1}+
\notag \\
B^{2}C_{y}^{1}+B^{2}C_{z}^{2}+B^{1}C_{y}^{2}-B^{1}C_{z}^{1}-2B^{2}B_{x}^{1}-2B^{1}B_{x}^{2}-2B_{xy}^{2}+2B_{xz}^{1}=0,
\tag{40}
\end{gather}%
\begin{gather}
6D^{1}A_{x}^{1}-6D^{2}A_{x}^{2}-3D^{1}B_{y}^{1}-3D^{1}B_{z}^{2}+3D^{2}B_{y}^{2}-3D^{2}B_{z}^{1}+3A^{1}D_{x}^{1}-3A^{2}D_{x}^{2}+
\notag \\
B_{xx}^{1}-2C_{xy}^{1}2C_{xz}^{2}-3B^{1}D_{y}^{1}-3B^{1}D_{z}^{2}+3B^{2}D_{z}^{2}-3B^{2}D_{z}^{1}+3D_{yy}^{1}-
\notag \\
3D_{zz}^{1}+6D_{yz}^{2}+2C^{1}C_{y}^{1}+2C^{1}C_{z}^{2}-2C^{2}C_{y}^{2}+2C^{2}C_{z}^{1}-C^{1}B_{x}^{1}+C^{2}B_{x}^{2}=0,
\tag{41}
\end{gather}%
\begin{gather}
6D^{2}A_{x}^{1}+6D^{1}A_{x}^{2}-3D^{2}B_{y}^{1}-3D^{2}B_{z}^{2}-3D^{1}B_{y}^{2}+3D^{1}B_{z}^{1}+3A^{2}D_{x}^{1}+3A^{1}D_{x}^{2}+
\notag \\
B_{xx}^{2}-2C_{xy}^{2}+2C_{xz}^{1}-3B^{2}D_{y}^{1}-3B^{2}D_{z}^{2}-3B^{1}D_{y}^{2}+3B^{1}D_{z}^{1}+3D_{yy}^{2}-
\notag \\
3D_{zz}^{2}-6D_{yz}^{1}+2C^{2}C_{y}^{1}-2C^{2}C_{z}^{2}+2C^{1}C_{y}^{2}-2C^{1}C_{z}^{1}-C^{2}B_{x}^{1}-C^{1}B_{x}^{2}=0;
\tag{42}
\end{gather}%
\textbf{3.} The system of RODEs corresponding to a r-CODE has four real
symmetries $\mathbf{X}_{1},~\mathbf{Y}_{1},~\mathbf{X}_{2}$ and $\mathbf{Y}%
_{2}$ with 
\begin{equation}
\mathbf{X}_{1}=\rho _{1}\mathbf{X}_{2}-\rho _{2}\mathbf{Y}_{2},~\mathbf{Y}%
_{1}\ =\rho _{1}\mathbf{Y}_{2}+\rho _{2}\mathbf{X}_{2},  \tag{43}
\end{equation}%
for nonconstant real functions $\rho _{1}(x)$ and $\rho _{2}(x)$ and they
satisfy%
\begin{equation}
\lbrack \mathbf{X}_{1},\mathbf{X}_{2}]-[\mathbf{Y}_{1},\mathbf{Y}_{2}]=0,~[%
\mathbf{X}_{1},\mathbf{Y}_{2}]+[\mathbf{Y}_{1},\mathbf{X}_{2}]=0,  \tag{44}
\end{equation}%
such that a point transformation $\chi =\chi (x,y,z),\Upsilon =\Upsilon
(x,y,z)$ and $\zeta =\zeta (x,y,z),~$which brings $\mathbf{X}_{1},~\mathbf{Y}%
_{1},~\mathbf{X}_{2}$ and $\mathbf{Y}_{2}$ to their canonical form%
\begin{equation}
\mathbf{X}_{1}=\frac{\partial }{\partial \Upsilon },~\mathbf{Y}_{1}=-\frac{%
\partial }{\partial \zeta },~\mathbf{X}_{2}=\chi \frac{\partial }{\partial
\Upsilon },~\mathbf{Y}_{2}=-\chi \frac{\partial }{\partial \zeta },  \tag{45}
\end{equation}%
reduces the system (35) and (36) to the linear form%
\begin{equation}
\Upsilon ^{\prime \prime }=W_{1}(\chi ),\text{ \ \ }\zeta ^{\prime \prime
}=W_{2}(\chi );  \tag{46}
\end{equation}%
\textbf{4.} The system of RODEs corresponding to a r-CODE has four real
symmetries $\mathbf{X}_{1},~\mathbf{Y}_{1},~\mathbf{X}_{2}$ and $\mathbf{Y}%
_{2}$ with 
\begin{equation}
\mathbf{X}_{1}=\rho _{1}\mathbf{X}_{2}-\rho _{2}\mathbf{Y}_{2},~\mathbf{Y}%
_{1}\ =\rho _{1}\mathbf{Y}_{2}+\rho _{2}\mathbf{X}_{2},  \tag{47}
\end{equation}%
for nonconstant $\rho _{1}$ and $\rho _{2}$ and they satisfy either 
\begin{equation}
\lbrack \mathbf{X}_{1},\mathbf{X}_{2}]-[\mathbf{Y}_{1},\mathbf{Y}_{2}]\neq 0%
\text{ or }[\mathbf{X}_{1},\mathbf{Y}_{2}]+[\mathbf{Y}_{1},\mathbf{X}%
_{2}]\neq 0,  \tag{48}
\end{equation}%
such that a point transformation $\chi =\chi (x,y,z),\Upsilon =\Upsilon
(x,y,z)$ and $\zeta =\zeta (x,y,z),~$which brings $\mathbf{X}_{1},~\mathbf{Y}%
_{1},~\mathbf{X}_{2}$ and $\mathbf{Y}_{2}$ to their canonical form%
\begin{equation}
\mathbf{X}_{1}=\frac{\partial }{\partial \Upsilon },~\mathbf{Y}_{1}=\frac{%
\partial }{\partial \zeta },~\mathbf{X}_{2}=\Upsilon \frac{\partial }{%
\partial \Upsilon }+\zeta \frac{\partial }{\partial \zeta },~\mathbf{Y}%
_{2}=\zeta \frac{\partial }{\partial \Upsilon }-\Upsilon \frac{\partial }{%
\partial \zeta },  \tag{49}
\end{equation}%
reduces the system (35) and (36) to the linear form%
\begin{equation}
\Upsilon ^{\prime \prime }=\Upsilon ^{\prime }W_{1}(\chi )-\zeta ^{\prime
}W_{2}(\chi ),  \tag{50}
\end{equation}%
\begin{equation}
\zeta ^{\prime \prime }=\Upsilon ^{\prime }W_{2}(\chi ).+\zeta ^{\prime
}W_{1}(\chi );  \tag{51}
\end{equation}%
\textbf{5.} The system of RODEs corresponding to a r-CODE has four real
symmetries $\mathbf{X}_{1},~\mathbf{X}_{2}$ and $\mathbf{Y}_{2}$ with 
\begin{equation}
\mathbf{X}_{1}\neq \rho _{1}\mathbf{X}_{2}-\rho _{2}\mathbf{Y}_{2},  \tag{52}
\end{equation}%
for nonconstant $\rho _{1}$ and $\rho _{2}$ and they satisfy%
\begin{equation}
\lbrack \mathbf{X}_{1},\mathbf{X}_{2}]=0,~[\mathbf{X}_{1},\mathbf{Y}_{2}]=0,
\tag{53}
\end{equation}%
such that a point transformation $\chi =\chi (x,y,z),\Upsilon =\Upsilon
(x,y,z)$ and $\zeta =\zeta (x,y,z),~$which brings $\mathbf{X}_{1},~\mathbf{X}%
_{2}$ and $\mathbf{Y}_{2}$ to their canonical form%
\begin{equation}
\mathbf{X}_{1}=\frac{\partial }{\partial \chi },\mathbf{X}_{2}=\frac{%
\partial }{\partial \Upsilon },~\mathbf{Y}_{2}=\frac{\partial }{\partial
\zeta },  \tag{54}
\end{equation}%
reduces the system (35) and (36) to the system of RODEs corresponding to
r-CODE is at most cubic in all its first derivatives;\newline
\textbf{6.} The system of RODEs corresponding to a r-CODE has four real
symmetries $\mathbf{X}_{1},~\mathbf{Y}_{1},~\mathbf{X}_{2}$ and $\mathbf{Y}%
_{2}$ with 
\begin{equation}
\mathbf{X}_{1}\neq \rho _{1}\mathbf{X}_{2}-\rho _{2}\mathbf{Y}_{2},~\mathbf{Y%
}_{1}\ \neq \rho _{1}\mathbf{Y}_{2}+\rho _{2}\mathbf{X}_{2},  \tag{55}
\end{equation}%
for nonconstant $\rho _{1}$ and $\rho _{2}$ and they satisfy either 
\begin{equation}
\lbrack \mathbf{X}_{1},\mathbf{X}_{2}]-[\mathbf{Y}_{1},\mathbf{Y}_{2}]\neq 0%
\text{ or }[\mathbf{X}_{1},\mathbf{Y}_{2}]+[\mathbf{Y}_{1},\mathbf{X}%
_{2}]\neq 0,  \tag{56}
\end{equation}%
such that a point transformation $\chi =\chi (x,y,z),\Upsilon =\Upsilon
(x,y,z)$ and $\zeta =\zeta (x,y,z)~,$which brings $\mathbf{X}_{1},~\mathbf{Y}%
_{1},~\mathbf{X}_{2}$ and $\mathbf{Y}_{2}$ to their canonical form%
\begin{equation}
\mathbf{X}_{1}=\frac{\partial }{\partial \Upsilon },~\mathbf{Y}_{1}=\frac{%
\partial }{\partial \zeta },~\mathbf{X}_{2}=\chi \frac{\partial }{\partial
\chi }+\Upsilon \frac{\partial }{\partial \Upsilon }+\zeta \frac{\partial }{%
\partial \zeta },{}  \tag{57}
\end{equation}%
\begin{equation}
~\mathbf{Y}_{2}=\zeta \frac{\partial }{\partial \Upsilon }-\Upsilon \frac{%
\partial }{\partial \zeta },  \tag{58}
\end{equation}%
reduces the system (35) and (36) to the linear form%
\begin{gather}
\chi \Upsilon ^{\prime \prime }=a_{1}(\Upsilon ^{\prime 3}-3\Upsilon
^{\prime }\zeta ^{\prime 2})-a_{2}(3\Upsilon ^{\prime 2}\zeta ^{\prime
}-\zeta ^{\prime 3})+b_{1}(\Upsilon ^{\prime 2}-\zeta ^{\prime
2})-2b_{2}\Upsilon ^{\prime }\zeta ^{\prime }+  \notag \\
\frac{1}{3(a_{1}^{2}+a_{2}^{2})}%
\{3(a_{1}^{2}+a_{2}^{2})+(b_{1}^{2}-b_{2}^{2})a_{1}+2b_{1}b_{2}a_{2}\}%
\Upsilon ^{\prime }-\frac{1}{3(a_{1}^{2}+a_{2}^{2})}\{2b_{1}b_{2}a_{1}- 
\notag \\
a_{2}(b_{1}^{2}-b_{2}^{2})\}\zeta ^{\prime }+\frac{1}{3(a_{1}^{2}+a_{2}^{2})}%
(b_{1}a_{1}+b_{2}a_{2})+\frac{1}{27(a_{1}^{2}+a_{2}^{2})}%
\{(b_{1}^{3}-3b_{1}b_{2}^{2})(a_{1}^{2}-a_{2}^{2})  \notag \\
+2a_{1}a_{2}(3b_{1}^{2}b_{2}-b_{2}^{3}),  \tag{59}
\end{gather}%
\begin{gather}
\chi \zeta ^{\prime \prime }=a_{2}(\Upsilon ^{\prime 3}-3\Upsilon ^{\prime
}\zeta ^{\prime 2})+a_{1}(3Y^{\prime 2}\zeta ^{\prime }-\zeta ^{\prime
3})+b_{2}(\Upsilon ^{\prime 2}-\zeta ^{\prime 2})+2b_{1}\Upsilon ^{\prime
}\zeta ^{\prime }+  \notag \\
\frac{1}{3(a_{1}^{2}+a_{2}^{2})}%
\{3(a_{1}^{2}+a_{2}^{2})+(b_{1}^{2}-b_{2}^{2})a_{1}+2b_{1}b_{2}a_{2}\}\zeta
^{\prime }+\frac{1}{3(a_{1}^{2}+a_{2}^{2})}\{2b_{1}b_{2}a_{1}-  \notag \\
a_{2}(b_{1}^{2}-b_{2}^{2})\}\Upsilon ^{\prime }+\frac{1}{%
3(a_{1}^{2}+a_{2}^{2})}(b_{2}a_{1}-b_{1}a_{2})+\frac{1}{%
27(a_{1}^{2}+a_{2}^{2})}\{(3b_{1}^{2}b_{2}-b_{2}^{3})(a_{1}^{2}-a_{2}^{2}) 
\notag \\
-2(b_{1}^{3}-3b_{1}b_{2}^{2})a_{1}a_{2}\}.  \tag{60}
\end{gather}%
The following invertible real transformation%
\begin{equation}
\tilde{\chi}=\Upsilon +\frac{1}{3(a_{1}^{2}+a_{2}^{2})}%
\{(b_{1}a_{1}+b_{2}a_{2})\chi \},  \tag{61}
\end{equation}%
\begin{equation}
\tilde{\Upsilon}=\frac{1}{2}(\Upsilon ^{2}-\zeta ^{2})+\frac{1}{%
3(a_{1}^{2}+a_{2}^{2})}[\{(b_{1}a_{1}+b_{2}a_{2})\chi \}\Upsilon
-\{(b_{2}a_{1}-b_{1}a_{2})\chi \}\zeta ]+  \notag
\end{equation}%
\begin{equation}
\frac{1}{18(a_{1}^{2}+a_{2}^{2})^{2}}%
[\{(b_{1}^{2}-b_{2}^{2})(a_{1}^{2}-a_{2}^{2})+4b_{1}b_{2}a_{1}a_{2}\}\chi
^{2}]+\frac{1}{2(a_{1}^{2}+a_{2}^{2})}\{a_{1}\chi ^{2}\},  \tag{62}
\end{equation}%
\begin{equation}
\tilde{\zeta}=\Upsilon \zeta +\frac{1}{3(a_{1}^{2}+a_{2}^{2})}%
[\{(b_{1}a_{1}+b_{2}a_{2})\chi \}\zeta +\{(b_{2}a_{1}-b_{1}a_{2})\chi \}F]+ 
\notag
\end{equation}%
\begin{equation}
\frac{1}{18(a_{1}^{2}+a_{2}^{2})^{2}}[-\chi
^{2}\{2b_{1}b_{2}(a_{1}^{2}-a_{2}^{2})-2a_{1}a_{2}(b_{1}^{2}-b_{2}^{2})\}]+%
\frac{1}{2(a_{1}^{2}+a_{2}^{2})}\{-a_{2}\chi ^{2}\},  \tag{63}
\end{equation}%
transforms (59) and (60) into the system of linear RODEs 
\begin{equation}
\Upsilon ^{\prime \prime }=0,\text{ \ \ }\zeta ^{\prime \prime }=0.  \tag{64}
\end{equation}%
\textbf{Examples.} Now we discuss some illustrative examples of linearizable
r-CODEs giving linearizable two dimensional systems of RODEs.\newline
\textbf{1.} Consider the second order nonlinear r-CODE with an arbitrary
function $w(x)$%
\begin{equation}
u~u^{\prime \prime }=u^{\prime 2}+w(x)u^{2}.  \tag{65}
\end{equation}%
The above r-CODE admits two CLSs of the form%
\begin{equation}
\mathbf{Z}_{1}=xu\frac{\partial }{\partial u},~\mathbf{Z}_{2}=u\frac{%
\partial }{\partial u},  \tag{66}
\end{equation}%
and thus%
\begin{equation}
\lbrack \mathbf{Z}_{1},\mathbf{Z}_{2}]=0\text{ and }\mathbf{Z}_{2}=\frac{1}{x%
}\mathbf{Z}_{1}.  \tag{67}
\end{equation}%
By using the complex transformation 
\begin{equation}
\chi =\frac{1}{x}\text{ and }U=\frac{1}{x}\log u,  \tag{68}
\end{equation}%
(65) can be reduced into a linear r-CODE%
\begin{equation}
U^{\prime \prime }=\frac{1}{\chi ^{3}}w(\frac{1}{\chi }).  \tag{69}
\end{equation}%
The system of RODEs corresponding to (65) are 
\begin{equation}
yy^{\prime \prime }-zz^{\prime \prime } =y^{\prime 2}-z^{\prime
2}+w_{1}(y^{2}-z^{2})-2yzw_{2},  \tag{70} \\
\end{equation}%
\begin{equation}
yz^{\prime \prime }+{\small z}y^{\prime \prime } =2y^{\prime }z^{\prime
}+2w_{1}yz+w_{2}(y^{2}-z^{2}),  \tag{71}
\end{equation}%
with coefficients that satisfy the the Lie type conditions $(39)-(42)$.
Invoking the transformation 
\begin{equation}
\chi =\frac{x}{x^{2}+y^{2}}\text{, \ }\Upsilon =\frac{1}{2x}\ln
(y^{2}+z^{2}),\text{ \ }\zeta =\frac{1}{x}\tan ^{-1}(\frac{z}{y}),  \tag{72}
\end{equation}%
equations (70) and (71) reduce to the linear system of RODEs 
\begin{equation}
\Upsilon ^{\prime \prime }=\frac{1}{\chi }w_{1},\text{ \ }\zeta ^{\prime
\prime }=\frac{1}{\chi }w_{2},  \tag{73}
\end{equation}%
where%
\begin{equation}
w_{1}=w_{1}(\frac{1}{\chi }),\text{ }w_{2}=w_{2}(\frac{1}{\chi }).  \tag{74}
\end{equation}%
\newline
\textbf{2.} Consider the nonlinear r-CODE 
\begin{equation}
u^{\prime \prime }+3u~u^{\prime }+u^{3}=0,  \tag{75}
\end{equation}%
which is linearizable as it satisfies the Lie-type conditions. The set of
two noncommuting CLSs are%
\begin{equation}
\mathbf{Z}_{1}=\frac{\partial }{\partial x},~\mathbf{Z}_{2}=x\frac{\partial 
}{\partial x}-u\frac{\partial }{\partial u}.  \tag{76}
\end{equation}%
Note that $\mathbf{X}_{1}\neq \rho (x,u)\mathbf{X}_{2}.$ Thus, we invoke
condition nine of Theorem 1 to find a linearization transformation. The
complex point transformation that reduces the symmetries (76) to their
canonical form is 
\begin{equation}
\chi =\frac{1}{u},~U=x+\frac{1}{u}  \tag{77}
\end{equation}%
and (75) reduces to 
\begin{equation}
\chi U^{\prime \prime }=-U^{\prime 3}+6U^{\prime 2}-11U^{\prime }+6  \tag{78}
\end{equation}%
by means of the transformation (77). Equation (75) linearizes to $\tilde{U}%
^{\prime \prime }=0,$ via the complex transformations with $a=-1,~b=6,$
and by using (33). That is 
\begin{equation}
\tilde{\chi}=x-\frac{1}{u},~\tilde{U}=\frac{x^{2}}{2}-\frac{x}{u}.  \tag{79}
\end{equation}%
The transformation (79) may seem strange as the variable $x$ is real while $\chi$ is 
complex. The point is that we are actually working with a complex independent
variable {\it restricted to lie on the real line}. We further require that
{\it at the end} the variable again be restricted to the real line, but in the 
intervening steps the variable moves off the real line. The procedure is reminiscent 
of analytic continuation. This odd behavior, mutually, does not appear in P - II.
In the new coordinates $(\tilde{\chi},\tilde{U})$, the solution of (75) is 
\begin{equation}
\tilde{U}=\alpha \tilde{\chi}+\beta ,  \tag{80}
\end{equation}%
where $\alpha $ and $\beta $ are complex constants whereas in coordinates $%
(\chi ,u)$ the above equation yields%
\begin{equation}
u=\frac{2(x-\alpha )}{x^{2}-2\alpha x-2\beta },  \tag{81}
\end{equation}%
which satisfies r-CODE (75). The system of RODEs corresponding to (75) is 
\begin{equation}
y^{\prime \prime }=-3(yy^{\prime }-zz^{\prime })-(y^{3}-3yz^{2}),  \tag{82}
\end{equation}%
\begin{equation}
z^{\prime \prime }=-3(zy^{\prime }+yz^{\prime })-(3y^{2}z-z^{3}),  \tag{83}
\end{equation}%
with coefficients that satisfy the Lie type conditions $(34)-(37),$ thus the
above system is linearizable. The general solution of system (82) and (83)
is 
\begin{equation}
y =\frac{2(x-\alpha _{1})(x^{2}-2\alpha _{1}x-2\beta _{1})+4\alpha
_{2}(\alpha _{2}x+\beta _{2})}{(x^{2}-2\alpha _{1}x-2\beta
_{1})^{2}+(2\alpha _{2}x+2\beta _{2})^{2}},  \tag{84} \\
\end{equation}%
\begin{equation}
z =\frac{4(x-\alpha _{1})(\alpha _{2}x+\beta _{2})-2\alpha
_{2}(x^{2}-2\alpha _{1}x-2\beta _{1})}{(x^{2}-2\alpha _{1}x-2\beta
_{1})^{2}+(2\alpha _{2}x+2\beta _{2})^{2}},  \tag{85}
\end{equation}%
which is obtained from (81). 
\newline
\newline
\textbf{3. }Consider the nonlinear r-CODE%
\begin{equation}
u^{\prime \prime }=1+(u^{\prime }-x)^{2}w(2u-x^{2})  \tag{86}
\end{equation}%
which admits the CLSs%
\begin{equation}
\mathbf{Z}_{1}=\frac{\partial }{\partial x}+x\frac{\partial }{\partial u},~%
\mathbf{Z}_{2}=x\frac{\partial }{\partial x}+x^{2}\frac{\partial }{\partial u%
}.  \tag{87}
\end{equation}%
Equation (86) reduces to a linear r-CODE by using the complex
transformation%
\begin{equation}
\chi =2u-x^{2}\text{ and }U=x,  \tag{88}
\end{equation}%
to become%
\begin{equation}
U^{\prime \prime }=\frac{-1}{2}U^{\prime }w(\chi ).  \tag{89}
\end{equation}%
Again, notice that (88) is not the usual complex
transformation as it is from (real, complex) to (complex, real) leading to
an apparent contradiction. As before, it corresponds to analytic continuation. To check
consistency, we may take $w=1$ i.e. $w_{1}=1,w_{2}=0,$ to obtain
\begin{equation}
u=\alpha +\ln 2+\frac{x^{2}}{2}-\log (\beta -x),  \tag{90}
\end{equation}%
where $\alpha $ and $\beta $ are complex constants. Note that (90) satisfy
(86). The system of RODEs corresponding to (86) is 
\begin{gather}
y^{\prime \prime }=1+\{(y^{\prime }-x)^{2}-z^{\prime
}{}^{2}\}w_{1}-2(y^{\prime }-x)z^{\prime }w_{2},  \tag{91} \\
z^{\prime \prime }=2(y^{\prime }-x)z^{\prime }w_{1}+\{(y^{\prime
}-x)^{2}-z^{\prime }{}^{2}\}w_{2},  \tag{92}
\end{gather}%
with coefficients satisfying the Lie type conditions $(39)-(42),$ where%
\begin{equation}
w_{1}=w_{1}(2y-x^{2},2z)\text{ and }w_{2}=w_{2}(2y-x^{2},2z).  \tag{93}
\end{equation}%
By placing values of $w_{1}=1, w_{2}=0,$ we get the system%
\begin{equation}
y^{\prime \prime }\text{ } =\text{ }1+\{(y^{\prime }-x)^{2}-z^{\prime
}{}^{2}\},  \tag{94} \\
\end{equation}%
\begin{equation}
z^{\prime \prime }\text{ } =\text{ }2(y^{\prime }-x)z^{\prime }, 
\tag{95}
\end{equation}%
with general solution 
\begin{equation}
y =\alpha _{1}-\ln 2+\frac{x^{2}}{2}-\frac{1}{2}\log ((\beta
_{1}-x)^{2}+\beta _{2}^{2}),  \tag{96} \\
\end{equation}%
\begin{equation}
z =\alpha _{2}-\tan ^{-1}(\frac{\beta _{2}}{\beta _{2}-x}).  \tag{97}
\end{equation}%
Also, if we take 
\begin{equation}
w(2u-x^{2})=\frac{1}{2u-x^{2}}, 
\tag{98}
\end{equation}%
the solution of (86) is 
\begin{equation}
u=\frac{x^{2}}{2}+\frac{1}{\sqrt{2}}\sqrt{\frac{\beta -x}{\alpha }}, 
\tag{99}
\end{equation}%
which gives the solution 
\begin{equation}
y =\frac{x^{2}}{2}+\frac{1}{\sqrt{2(\alpha _{1}^{2}+\alpha _{2}^{2})}}%
\sqrt{(\alpha _{1}(\beta _{1}-x)+\beta _{2}\alpha _{2})^{2}+(\beta
_{2}\alpha _{2}-\alpha _{2}(\beta _{1}-x))^{2}}\times   \notag \\
\end{equation}%
\begin{equation}
\cos (\frac{\beta _{2}\alpha _{2}-\alpha _{2}(\beta _{1}-x)}{2(\alpha
_{1}(\beta _{1}-x)+\beta _{2}\alpha _{2})}),  \tag{100} \\
\end{equation}%
\begin{equation}
z =\frac{1}{\sqrt{2(\alpha _{1}^{2}+\alpha _{2}^{2})}}\sqrt{(\alpha
_{1}(\beta _{1}-x)+\beta _{2}\alpha _{2})^{2}+(\beta _{2}\alpha _{2}-\alpha
_{2}(\beta _{1}-x))^{2}}\times   \notag \\
\end{equation}%
\begin{equation}
\sin (\frac{\beta _{2}\alpha _{2}-\alpha _{2}(\beta _{1}-x)}{2(\alpha
_{1}(\beta _{1}-x)+\beta _{2}\alpha _{2})}),  \tag{101}
\end{equation}%
of system 
\begin{equation}
y^{\prime \prime } {\small =}1+\{(y^{\prime }-x)^{2}-z^{\prime }{}^{2}\}%
\frac{2y-x^{2}}{(2y-x^{2})^{2}+4z^{2}}+4(y^{\prime }-x)z^{\prime }\frac{z}{%
(2y-x^{2})^{2}+4z^{2}},  \tag{102} \\
\end{equation}%
\begin{equation}
z^{\prime \prime } {\small =}2(y^{\prime }-x)z^{\prime }\frac{2y-x^{2}}{%
(2y-x^{2})^{2}+4z^{2}}-2\{(y^{\prime }-x)^{2}-z^{\prime }{}^{2}\}\frac{z}{%
(2y-x^{2})^{2}+4z^{2}}.  \tag{103}
\end{equation}%
It is interesting to see that here $w$ is an arbitrary
complex function that give a class of systems of RODEs that
correspond to r-CODE (86). Thus, the linearization of a general r-CODE like
(86) encodes the linearization of a large class of systems of RODEs.
Further, the complex solutions like (90) directly give us the real solutions
of those systems of RODEs yielding non-trivial examples of complex magic.

\section{\textbf{Summary and Discussion}}

\qquad The study of linearization of RODEs is significant as it plays an
important role in the reduction of equations to simple form, from which one
can construct exact solutions. It also includes finding the transformations
that map nonlinear RODEs into linear RODEs. Linearization of systems of
nonlinear RODEs is far from trivial as the Lie conditions are very
complicated and involve a system of more than ten equations to be solved
simultaneously $[14-17]$.

Our approach helps in constructing non-trivial ways of \textit{reduction in
order, linearization} \textit{and solving }of certain systems of RODEs. The
CLS analysis gets used elegantly and several known results in classical
symmetry analysis are extended to complex functions of real variables in
this paper. The idea of CLSs and their use in CODEs is discussed in \cite%
{Ali}. It is proved that one can use CLSs for studying RLSs for systems of
PDEs corresponding to CODEs. We have also studied the use of CLSs in
variational problems in solving inverse problems which is to find a complex
Lagrangian for CODEs which results in real Lagrangians of systems of
corresponding PDEs and double reduction via complex Noether symmetries in 
\cite{Ali1}. The CLS analysis is further applied to solve the inverse
problems for systems of ODE in \cite{Ali3}

In this work we have looked at linearizability criteria for second order
r-CODEs. We obtained analogous Lie type conditions for r-CODEs and applied
them for linearization of system of two nonlinear second order RODEs. The
corresponding statement for a system of two RODEs associated to a second
order r-CODE was presented. Examples were given for the linearization of
systems of RODEs including the construction of the point transformations
that did the reduction to linear form. The procedure adopted, of taking a
linearizable RODE, by converting an r-CODE and thence obtaining a system of
two linearizable RODEs, has provided a powerful and elegant alternative way
of analyzing systems of RODEs. Further, the complex solutions of r-CODEs
directly give the real solutions of systems of RODEs corresponding to them.

An extensive classification of those systems of RODEs that correspond to
r-CODEs would be of great worth. It is important to point out the fact that
the system $(34)-(37)$ provides a necessary and sufficient condition for
system (32) and (33) to be linearizable. Lie conditions for a general system
of second order quadratically semi-linear ordinary differential equations
were obtained in \cite{mah}. Whether system (32) and (33) constitute a
special case of a system of second order cubically semi-linear RODEs in \cite%
{mah2} and whether conditions $(34)-(37)$ can be derived from conditions in 
\cite{mah2} remains to be explored. Whether the solutions obtained in \cite%
{mah2} for a system of second order cubically semilinear RODEs can be
derived from complex solutions and then what is the difference in the real
transformations obtained from complex transformations from those real
transfomation mentioned in \cite{mah2} for linearization are to be answered.
It may happen that in some cases one approach becomes a special case of the
other and vice versa in other cases.

An important feature of complex magic in r-CODEs that remains to be explored completely 
is the procedure that appears to be analogous to analytic continuation. Explicitly,
we dealt with complex functions of one real variable that yield two real
functions of one real variable, although, the transformation that linearizes the
r-CODEs (70) and (84) was not the usual complex transfomation. Thus, in
a way we are mapping real solutions into real solutions via complex mapping.
Due to complex integration of r-CODEs the constants in the solutions are
also complex which yields two real constants when complex solutions are
transformed into real solutions. In this paper, the complex function is
restricted to depend on the real line instead of lying on the entire complex plane.
This procedure does not arise in P - II.

It is suggested that Lie criteria for linearization of two r-CODEs can be
found in a similar way, which can then be used to linearize four nonlinear
RODEs corresponding to a system of two r-CODEs. It is well known that a
system of two second order RODES admits $5,6,7,8$ or $15$ RLSs and the
maximal symmetry algebra is $sl(4,R)$ for the simplest system \cite{waf}.
Thus, it is conjectured that the system of two second order r-CODEs admits $%
10,12,14,16$ or $30$ RLSs corresponding to $5,6,7,8$ or $15$ CLSs. The RLSs $%
10,12,14,16$ or $30$ are basically the symmetries of four second order RODEs
that correspond to two second order r-CODEs. Further, the simplest systems
of two second order r-CODEs admits a maximal symmetry algebra $sl(4,C)$. The
earlier procedures $[14-17]$ already became unwieldy for systems of four
RODEs. It seems likely that the system of two r-CODEs may provide a
convenient way to obtain linearization of systems of four RODEs. Further,
the extraction of real solutions of systems of four RODE from complex
solutions of system of two r-CODEs is required to be analyzed.

One can construct Lie type criteria for a third order r-CODE, which will
help in the linearization of systems of two nonlinear RODEs of third order.
Previous results for third and fourth order RODE\ \cite{mah3, mah5} might be
extendible by this method to systems of RODEs and PDEs. It is hoped that CLS
analysis would provide linearizability conditions similar to (34)-(37) in a
direct way for systems of RODEs instead of larger number of linearizability
conditions mentioned in \cite{mah3, mah5} that makes algebraic calculations
very difficult and tedious. Also, classification of systems of third order
RODEs that correspond to r-CODEs and their systems would be worth exploring.
Further, it is interesting to analyze how the complex solutions of third
order r-CODEs would yield real solutions of systems of third order RODEs.

Another application to obtain the generalization of Lie type conditions for
PDEs, is investigated in P - II \cite{Ali2}.

\section{\textbf{Acknowledgments}}

\qquad SA is most grateful to NUST and DECMA in providing finnancial
assistance for his stay in Wits University, South Africa, where this work
was done. SA would also like to acknowledge Roger Penrose's book `\textit{%
The Road to Reality}' for making him believe in the nontrivial implications
of the magic in complex theory. AQ is grateful to the School of
Computational and Applied Mathematics (DECMA) for funding his stay in the
university.

\end{document}